%% file: main.tex
\documentclass[]{elsarticle}

\usepackage{subcaption,wrapfig}

\usepackage{amsthm,amsmath,amssymb,tikz,stmaryrd}
\usepackage{tikz}
\usetikzlibrary{fit, arrows, shapes, positioning, shadows, matrix, calc}
 
\usetikzlibrary{decorations.pathreplacing}
\tikzstyle std=[line width=0.7pt]   
\tikzstyle stdthin=[line width=0.3pt]   
\tikzstyle stdthick=[line width=1.0pt]   

\tikzstyle fwd=[line width=0.7pt, ->]   
\tikzstyle fwdthin=[line width=0.3pt, ->]   
\tikzstyle fwdthick=[line width=1.0pt, ->]   
\tikzstyle fwddash=[line width=0.7pt, dashed, ->]   

\tikzstyle bwd=[double, line width=0.3pt, ->]  
\tikzstyle refl=[double,  dashed, line width=0.2pt, ->]    

\tikzstyle{every node}=[font=\small] 

\tikzset{
  >=stealth', 
  invisible/.style={opacity=0}, 
  alt/.code args={<#1>#2#3}{\alt<#1>{\pgfkeysalso{#2}}{\pgfkeysalso{#3}}}, 
  visible on/.style={alt=#1{}{invisible}}, 
  smallnode/.style={circle, fill=black, thick, inner sep=1pt, minimum size=1.5pt}, 
  punkt/.style={
           rectangle,
           rounded corners,
           draw=black, very thick,
           text width=5.5em,
           minimum height=2em,
           text centered},
  punkt_big/.style={
           rectangle,
           rounded corners,
           draw=black, very thick,
           text width=7em,
           minimum height=2em,
           text centered},
}

  \usetikzlibrary{positioning}
  \usetikzlibrary{shadows}
  \usetikzlibrary{fit, arrows, shapes, positioning, shadows, matrix, calc}
\usetikzlibrary{decorations.pathreplacing}

\usetikzlibrary{shapes,calc,positioning}
\usetikzlibrary{calc,intersections}
\usepackage{tkz-euclide}
\usetkzobj{all}

\newcommand{\reals}{\mathbb{R}}
\newcommand{\complex}{\mathbb{C}}

\newcommand{\hilbert}{\mathcal{H}}

\input{FancyMacros}

\definecolor{lightgrey}{gray}{0.8}
\definecolor{medgrey}{gray}{0.6}
\definecolor{darkgrey}{gray}{0.4}

  \usetikzlibrary{positioning}
  \usetikzlibrary{shadows}
  \usetikzlibrary{fit, arrows, shapes, positioning, shadows, matrix, calc}
\usetikzlibrary{decorations.pathreplacing}

\usepackage{graphicx}

\newtheorem{theorem}{Theorem}

\newtheorem{corollary}{Corollary}
\newtheorem{fact}{Fact}

\newcommand{\C}{\mathrm{Circ}}

\newcommand{\fix}{\mathrm{Fix}}

\newcommand{\Avc}{\mathrm{Disk}}

\begin{document}

\begin{frontmatter}

\title{Tight Coefficients of Averaged Operators\\ via Scaled Relative Graph}

%
%
%

\author{Xinmeng Huang\fnref{xinmengaddress}}
\ead{hxm8888@mail.ustc.edu.cn}

\author{Ernest K. Ryu\fnref{ryuaddress}\corref{mycorrespondingauthor}}
\ead{ernestryu@snu.ac.kr}

\author{Wotao Yin\fnref{yinaddress}}
\ead{wotaoyin@math.ucla.edu}

\cortext[mycorrespondingauthor]{Corresponding author}

\address[xinmengaddress]{School of Mathematical Sciences, University of Science and Technology of China}
\address[ryuaddress]{Department of Mathematical Sciences, Seoul National University}
\address[yinaddress]{Department of Mathematics, University of California, Los Angeles}

%
%
%

\begin{abstract}
Many iterative methods in optimization are fixed-point iterations with averaged operators. As such methods converge at an $\mathcal{O}(1/k)$ rate with the constant determined by the averagedness coefficient, establishing small averagedness coefficients for operators is of broad interest. In this paper, we show that the averagedness coefficients of the composition of averaged operators by Ogura and Yamada (Numer Func Anal Opt 32(1--2):113--137, 2002)  and the three-operator splitting by Davis and Yin (Set-Valued Var Anal 25(4):829--858, 2017) are tight. The analysis relies on the scaled relative graph, a geometric tool recently proposed by Ryu, Hannah, and Yin (arXiv:1902.09788, 2019).
\end{abstract}


\end{frontmatter}


\section{Introduction}


%
Since their introduction in \cite{bailion1978asymptotic}, averaged operators have been widely used in the analysis of nonlinear fixed-point iterations.
In particular, a wide range of optimization methods can be analyzed as fixed-point iterations with a composition of averaged operators, which are themselves averaged \cite{combettes2004solving}.
The smallest (best) averagedness coefficient for this setup was presented by Ogura and Yamada \cite{ogura2002non} and was introduced to the broader optimization community by Combettes and Yamada \cite{combettes2015} and Bauschke and Combettes \cite{BCBook}.
More recently, Davis and Yin presented a three-operator splitting method and established its convergence by showing the associated operator is averaged \cite{davis2017three}.
Whether these averagedness coefficients are tight, loosely defined as being unable to be improved without additional assumptions, was not known.

The \emph{Scaled Relative Graph} (SRG) is a geometric tool for analyzing fixed-point iterations recently proposed by Ryu, Hannah, and Yin \cite{ryu2019scaled}.
The SRG maps the action of a nonlinear operator to a subset the 2D plane, analogous to how the spectrum maps the action of a linear operator to the complex plane.
A strength of the SRG is that it is well-suited for \emph{tight} analysis.

2D geometric illustrations have been used 
by Eckstein and Bertsekas \cite{eckstein1989,eckstein1992}, 
Giselsson \cite{giselsson2017linear,giselsson_slides},
Banjac and Goulart \cite{banjac2018},
and Giselsson and Moursi \cite{giselsson_moursi2019}
to qualitatively understand convergence of optimization algorithms.
The SRG is a rigorous formulation of such illustrations.



In this paper, we use the SRG to show tightness of the averagedness coefficients of the composition of averaged operators by Ogura and Yamada and the three-operator splitting by Davis and Yin.
Section~\ref{s:prelim} discusses general preliminaries and sets up the notation. Section~\ref{s:twoave} presents results on the composition of averaged operators. Section~\ref{s:dys} presents results on the Davis--Yin splitting.

\subsection{Contribution and prior work}
The contribution of this paper is in the results Corollaries~\ref{cor:ogura} and \ref{cor:dys}, which establish tightness of the averagedness coefficients, and the geometric proof technique based on the SRG.

The geometric arguments of Section~\ref{s:dys} are entirely new.
The geometric arguments of Section~\ref{s:twoave} overlap with the classical work on ``circular arithmetic'' initiated by Gargantini and Henrici \cite{Gargantini1971}.
In \cite{Hauenschild1974,hauenchild1980}, Hauenchild introduced the notion of ``optimal circular multiplication'', which considers the smallest circle enclosing the Minkowski product (defined in Section~\ref{s:prelim}) of two disks on the complex plane. This is \emph{not} the same as what we consider in Section~\ref{s:twoave}, since we find the smallest circle under the additional requirement that it goes through the point $(1,0)$. These two notions coincide sometimes, but not always. In \cite{polyak1994}, Polyak, Scherbakov, and Schmulyian perform calculations similar to that of Theorem~\ref{2-ave_cur} in the context of control theoretic stability analysis. In fact, Theorem 1 of \cite{polyak1994} is, after a change of variables, the same as Theorem~\ref{2-ave_cur} of this work. However, the proof in \cite{polyak1994} is not rigorous as it omits what we call Step 2 and Step 3 in our proof of Theorem~\ref{2-ave_cur}. In \cite{Farouki2001,farouki2002exact}, Farouki et al.\ also perform similar envelope calculations that are, after a change of variables, the same as that of Theorem~\ref{2-ave_cur} of this work. However, Farouki et al.\ also do not prove Steps 2 and 3; they merely state, without providing or outlining a proof, in Section 6.7 of \cite{Farouki2001} that ``one can easily see'' this fact.
To summarize, in the proof of Theorem~\ref{2-ave_cur}, Step 1 coincides with existing work, while Steps 2 and 3 are new.
Furthermore, the proof of Corollary~\ref{cor:ogura}, which connects the geometric analysis to the composition of averaged operators using the SRG, is new.





\section{Preliminaries}
\label{s:prelim}
We follow the standard notation of \cite{BCBook,ryu2016primer}. 
Write $\hilbert$ for a real Hilbert space equipped with the inner product $\langle\cdot,\cdot\rangle$ and norm $\|\cdot\|$. 
Write $A:\hilbert\rightrightarrows\hilbert$ to denote $A$ is a multi-valued operator on $\hilbert $. 
Write $I:\hilbert\rightarrow\hilbert$ for the identity operator.
Define the resolvent of $A$ as $J_A=(I+A)^{-1}$.
We say $A\colon\hilbert\rightrightarrows\hilbert$ is monotone if
\[
\langle u-v,x-y\rangle \ge 0,\qquad \forall x,y\in\hilbert, \,u\in Ax,\,v\in Ay.
\]
Write $\cM$ for the class of monotone operators.
For $\beta\in(0,\infty)$, we say a single-valued operator $A\colon\hilbert\to\hilbert$ is $\beta$-cocoercive if
\[
\langle Ax-Ay,x-y\rangle\ge \beta\|Ax-Ay\|^2,
\]
for all $x,y\in \hilbert$ and write $\cC_\beta$ for the class of $\beta$-cocoercive operators.
For $\theta\in(0,1)$, we say an operator $A$ is $\theta$-averaged if $A=(1-\theta)I+\theta N$ for some nonexpansive operator $N$ and write $\cN_\theta$ for the class of $\theta$-averaged operators.

We write complex numbers with the Cartesian and polar coordinate representations $z=x+yi$ and $z = re^{i \varphi}=r \cos(\varphi) + i r \sin(\varphi)$.
For notational convenience, we often identify $\complex$ with $\reals^2$.
We use Minkowski-type set notation with sets of complex numbers.
In particular, given $\alpha\in \complex$ and $Z,W\subseteq \complex$, write
\[
\alpha Z = \{ \alpha z\,|\,z\in Z\},\quad
ZW=\{zw\,|\,z\in Z,\,w\in W\}.
\]
The set $ZW$ is called the \emph{Minkowski product} of $Z$ and $W$.
Given a set $U$, write $\partial U$ to denote its boundary and write $U^\circ=U\backslash \partial U$ to denotes its interior.

\subsection{Scaled Relative Graph}
We follow the notation of \cite{ryu2019scaled}.
The scaled relative graph (SRG) of an operator $A$ is defined as 
\begin{align*}
\mathcal{G}(A)&=
\left\{
\frac{\|u-v\|}{\|x-y\|}
\exp\left[\pm i \angle (u-v,x-y)\right]
\,\Big|\,
u\in Ax,\,v\in Ay,\, x\ne y \right\}\\
&\qquad\qquad\qquad\qquad\qquad\qquad\qquad\qquad\qquad
\bigg(\cup \{\infty\}\text{ if $A$ is multi-valued}\bigg),
\end{align*}
where 
\begin{align*}
    \angle(a,b)=
    \left\{
    \begin{array}{ll}
    \arccos\p{\tfrac{\dotp{a,b}}{\n{a}\n{b}}}&\text{ if }a\ne 0,\,b\ne 0\\
    0&\text{ otherwise,}
    \end{array}
    \right.
\end{align*}
denotes the angle between $a,b\in\hilbert$.
The SRG $\cG(A)$ maps the action of the operator $A$ to points onto the extended complex plane.
The magnitude of each element of $\mathcal{G}(A)$, $\frac{\|u-v\|}{\|x-y\|}$, represents the size of the change in outputs $u,v$ relative to the size of the change in inputs $x,y$. The angle, $\angle (u-v,x-y)$, represents how much the change in outputs is aligned with the change in inputs.
The SRG of the class of operators $\cA$ is defined as
\[
\mathcal{G}(\cA)=\bigcup_{A\in \cA}\mathcal{G}(A).
\]


For $\theta\in(0,1)$, define
\[
\Avc(\theta)=\{z\in\complex\,\big|\,|z-(1-\theta)|\leq\theta\},
\quad\C(\theta)=\{z\in\complex\,\big|\,|z-(1-\theta)|=\theta\}.
\]
The sets have $(1-\theta)$ as their center and include $1$ as the right-most point.

\begin{fact}[Proposition~3.3 of \cite{ryu2019scaled}]
Let $\theta\in(0,1)$. Then
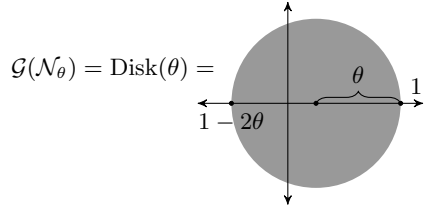
\begin{center}
\begin{tabular}{c}
\raisebox{-.5\height}{
\begin{tikzpicture}[scale=1.5]
\draw (-1.55,0.3) node {$\cG(\cN_\theta)=\Avc(\theta)=$};
\fill[fill=medgrey] (0.25,0) circle (0.75);
\draw [<->] (-.8,0) -- (1.2,0);
\draw [<->] (0,-0.9) -- (0,0.9);
\draw (1,-0) node [above right] {$1$};
\draw (-0.5,-0) node [below] {$1-2\theta$};
\filldraw (1,0) circle ({0.5*1.5/1.5pt});
\filldraw (.25,0) circle ({0.5*1.5/1.5pt});
\filldraw (-.5,0) circle ({0.5*1.5/1.5pt});
\draw (0.625,0.1) node [above] {$\theta$};
\draw [decorate,decoration={brace,amplitude=4.5pt}] (0.25,0) -- (1,0) ;
\end{tikzpicture}}
\end{tabular}
\end{center}
\end{fact}

\begin{fact}[Theorem~3.5 of \cite{ryu2019scaled}]
\label{fact:srg}
For the operator class $\cN_\theta$, where $\theta\in(0,1)$, inclusion within the operator class is equivalent to the inclusion of the SRG in the 2D plane.
\end{fact}

\subsection{Osculating circle, curvature, and envelope}
\label{sec_geom}
In differential geometry of curves, the \emph{osculating circle} of a sufficiently smooth plane curve $C$ at a point $P$ on the curve is the circle passing through $P$ that approximates $C$ most tightly within infinitesimal neighborhoods of $P$.
The center of the circle lies on the inner normal line, and the reciprocal of its radius is the \emph{curvature} of $C$ at $P$ \cite{coolidge1952unsatisfactory}.
For curves defined through polar coordinates as $r(\varphi)$, the curvature $\kappa(\varphi)$ at $r(\varphi)$ is given by \cite{pressley2010elementary}:
\begin{equation}
\kappa(\varphi)=\frac{r(\varphi)^2+2(\frac{dr}{d\varphi})^2-r(\varphi)\frac{d^2r}{d\varphi^2}}{(r(\varphi)^2+(\frac{dr}{d\varphi})^2)^{\frac{3}{2}}}
\label{eq:curvature}
\end{equation}
The osculating circle of $C$ at $P$ provides insight on the smallest circle through $P$ enclosing $C$.
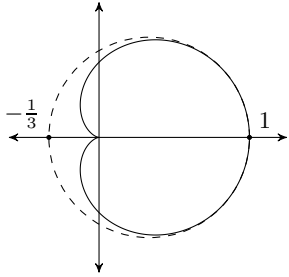
\begin{center}
\begin{tabular}{c}
\raisebox{-.5\height}{
\begin{tikzpicture}[scale=2]
\def\a{pi*.17};
\def\b{pi*.45};
\def\c{pi*.66};
\def\w{0.9};
\def\u{1.25};
\draw[dashed] (1/3,0) circle (2/3);
\draw[samples=200,smooth] plot[domain=-pi:pi] (xy polar cs:angle=\x r,radius= {cos(\x r/2)^2});
\filldraw (-1/3,0) circle[radius={0.5*1.5/2pt}];
\draw (-1/3,0) node [above left] {$-\frac{1}{3}$};
\filldraw (1,0) circle[radius={0.5*1.5/2pt}];
\draw (1,0) node [above right] {$1$};
\draw [<->] (-0.6,0) -- (\u,0);
\draw [<->] (0,-\w) -- (0,\w);
\draw (0.5,-1.1) node {$C=\left\{re^{i\varphi}\,|\,r= \frac{1}{2}(1+\cos(\varphi))\right\}$ and its osculating circle at 1};
\end{tikzpicture}}
\end{tabular}
\end{center}



An \emph{envelope} of a family of curves in the plane is a curve that is tangent to each member of the family at some point.
Formally, let each $\{C_t\}_{t\in\reals}$ be a parameterized family of curves in $\reals^2$ defined by $F(t,\mathbf{x})=0$, where $t\in\reals$ is the parameter, $\mathbf{x}\in\reals^2$, and $F$ is smooth. That is, $C_t = \{\mathbf{x}\,|\,F(t,\mathbf{x})=0\}$.
The envelope of $\{C_t\}_{t\in\reals}$ is defined as the set points satisfying 
\begin{equation}
F(t,\mathbf{x})=0,\qquad\frac{\partial F}{\partial t}(t,\mathbf{x})=0
\label{eq:envelope_formula}
\end{equation}
for some $t$.
The envelope includes the boundary of the region filled by the curves \cite[Section 5.17]{bruce1992curves}.
See \cite{boltyanskii1964envelopes,bruce1992curves} for further discussion.

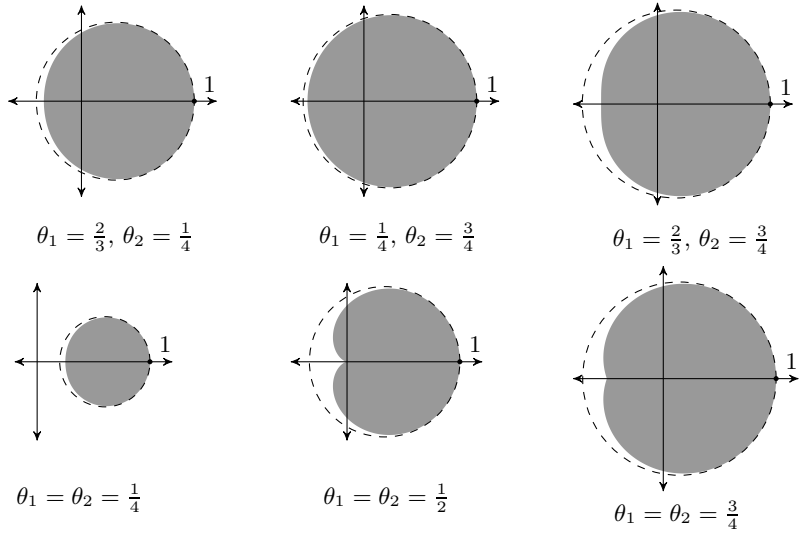
\begin{figure}
\centering
\begin{subfigure}[b]{0.3\textwidth}
\raisebox{-.5\height}{
\begin{tikzpicture}[scale=1.5]
\def\a{2/3};
\def\b{1/4};
\def\s{(\a+\b-2*\a*\b)/(1-\a*\b)};
\fill[fill=medgrey,samples=200,smooth] plot[domain=-180:180](\x:{cos(\x)*(1-\a)*(1-\b)+\a*\b+sqrt((cos(\x)*(1-\a)*(1-\b)+\a*\b)^2-(1-2*\a)*(1-2*\b))});
\draw[dashed] ({1-\s},0) circle ({\s});
\draw [<->] (-.65,0) -- (1.2,0);
\draw [<->] (0,-.85) -- (0,.85);
\draw (1,0)  node [above right] {1};
\filldraw (1,0) circle ({0.5*1.5/1.5pt});
\draw (0.3,-1) node [below] {$\theta_1=\frac{2}{3},\,\theta_2=\frac{1}{4}$};
\end{tikzpicture}
}
\end{subfigure}
\begin{subfigure}[b]{0.3\textwidth}
\raisebox{-.5\height}{
\begin{tikzpicture}[scale=1.5]
\def\a{1/4};
\def\b{3/4};
\def\s{(\a+\b-2*\a*\b)/(1-\a*\b)};
\fill[fill=medgrey,samples=200,smooth] plot[domain=-180:180](\x:{cos(\x)*(1-\a)*(1-\b)+\a*\b+sqrt((cos(\x)*(1-\a)*(1-\b)+\a*\b)^2-(1-2*\a)*(1-2*\b))});
\draw[dashed] ({1-\s},0) circle ({\s});
\draw [<->] (-.65,0) -- (1.2,0);
\draw [<->] (0,-.85) -- (0,.85);
\draw (1,0)  node [above right] {1};
\filldraw (1,0) circle ({0.5*1.5/1.5pt});
\draw (0.3,-1) node [below] {$\theta_1=\frac{1}{4},\,\theta_2=\frac{3}{4}$};
\end{tikzpicture}
}
\end{subfigure}
\begin{subfigure}[b]{0.3\textwidth}
\raisebox{-.5\height}{
\begin{tikzpicture}[scale=1.5]
\def\a{2/3};
\def\b{3/4};
\def\s{(\a+\b-2*\a*\b)/(1-\a*\b)};
\fill[fill=medgrey,samples=200,smooth] plot[domain=-180:180](\x:{cos(\x)*(1-\a)*(1-\b)+\a*\b+sqrt((cos(\x)*(1-\a)*(1-\b)+\a*\b)^2-(1-2*\a)*(1-2*\b))});
\draw[dashed] ({1-\s},0) circle ({\s});
\draw [<->] (-.75,0) -- (1.2,0);
\draw [<->] (0,-.9) -- (0,.9);
\draw (1,0)  node [above right] {1};
\filldraw (1,0) circle ({0.5*1.5/1.5pt});
\draw (0.3,-1) node [below] {$\theta_1=\frac{2}{3},\,\theta_2=\frac{3}{4}$};
\end{tikzpicture}
}
\end{subfigure}
\begin{subfigure}[b]{0.3\textwidth}
\raisebox{-.5\height}{
\begin{tikzpicture}[scale=1.5]
\def\t{0.25};
\def\s{2/5};
\fill[fill=medgrey,samples=200,smooth] plot[domain=-180:180](\x:{2*\t*(1-\t)+2*\t^2*cos(\x)});
\draw[dashed] ({2*\t-\s},0) circle (\s);
\draw [<->] (-0.2-1+2*\t,0) -- (1.2-1+2*\t,0);
\draw [<->] (-1+2*\t,-.7) -- (-1+2*\t,.7);
\draw (2*\t,0)  node [above right] {1};
\filldraw (2*\t,0) circle ({0.5*1.5/1.5pt});
\draw (2*\t^2-3*\t+0.5,-1) node [below] {$\theta_1=\theta_2=\frac{1}{4}$};
\end{tikzpicture}
}
\end{subfigure}
\begin{subfigure}[b]{0.3\textwidth}
\raisebox{-.5\height}{
\begin{tikzpicture}[scale=1.5]
\def\s{1/1.5};
\fill[fill=medgrey,samples=200,smooth] plot[domain=-180:180](\x:{1/2+1/2*cos(\x)});
\draw[dashed] ({1-\s},0) circle ({\s});
\draw [<->] (-.5,0) -- (1.2,0);
\draw [<->] (0,-.7) -- (0,.7);
\draw (1,0)  node [above right] {1};
\filldraw (1,0) circle ({0.5*1.5/1.5pt});
\draw (0.35,-1) node [below] {$\theta_1=\theta_2=\frac{1}{2}$};
\end{tikzpicture}
}
\end{subfigure}
\begin{subfigure}[b]{0.3\textwidth}
\raisebox{-.5\height}{
\begin{tikzpicture}[scale=1.5]
\def\t{0.75};
\def\s{(2*\t)/(1+\t)};
\fill[fill=medgrey,samples=200,smooth] plot[domain=-180:180](\x:{2*\t*(1-\t)+2*\t^2*cos(\x)});
\draw[dashed] ({2*\t-\s},0) circle ({\s});
\draw [<->] (-0.8-1+2*\t,0) -- (1.2-1+2*\t,0);
\draw [<->] (-1+2*\t,-1) -- (-1+2*\t,1);
\draw (2*\t,0)  node [above right] {1};
\filldraw (2*\t,0) circle ({0.5*1.5/1.5pt});
\draw (2*\t^2-2*\t+1,-1) node [below] {$\theta_1=\theta_2=\frac{3}{4}$};
\end{tikzpicture}
}
\end{subfigure}
\caption{
The shaded regions illustrate $\cG(\cN_{\theta_1}\cN_{\theta_2})$ given by Theorem~\ref{2-ave_cur}.
The circles drawn in dashed lines illustrate $\cG(\cN_\theta)$ given by Corollary~\ref{cor:ogura}.
}
\label{fig:srg-exmaples}
\end{figure}

\section{Tight characterization of the composition of averaged operators}
\label{s:twoave}
The composition of two averaged operators is itself an averaged operator, and Ogura and Yamada \cite{ogura2002non} showed the best known averagedness coefficient for this setup. In this section, we provide an alternate geometric proof of this result and establish its tightness.


Again, $\cN_{\theta_1}$ and $\cN_{\theta_2}$ are the classes of $\theta_1$- and $\theta_2$-averaged operators.
Define
$\cN_{\theta_1}\cN_{\theta_2}=\{N_1N_2\,|\,N_1\in \cN_{\theta_1},\,N_2\in \cN_{\theta_2}\}$
to be the class of compositions of $\theta_1$- and $\theta_2$-averaged operators.

\begin{theorem}
\label{2-ave_cur}
Let $\theta_1,\theta_2\in(0,1)$. 
Then $\mathcal{G}(\cN_{\theta_1}\cN_{\theta_2})$ is the region enclosed by the outer curve defined by
\begin{equation}
r(\varphi)^2-2r(\varphi)(\cos(\varphi)(1-\theta_1)(1-\theta_2)+\theta_1\theta_2)+(1-2\theta_1)(1-2\theta_2)=0.
\label{eq:srg-curve}
\end{equation}
\end{theorem}

Figure~\ref{fig:srg-exmaples} illustrates Theorem~\ref{2-ave_cur}.
To clarify, the equation of Theorem~\ref{2-ave_cur} defines at most two non-intersecting closed curves, one enclosing the other, and the SRG is given by the outer curve.

\begin{corollary}
\label{cor:ogura}
Let $\theta_1,\theta_2\in(0,1)$. Then $\cN_{\theta_1}\cN_{\theta_2}\subseteq \cN_{\theta}$ with
\[
\theta=\frac{\theta_1+\theta_2-2\theta_1\theta_2}{1-\theta_1\theta_2}.
\]
The averagedness coefficient $\theta$ is tight in the sense that it cannot be reduced without further assumptions.
\end{corollary}

\subsection{Proof of Theorem~\ref{2-ave_cur}}
Define $S$ to be the outer curve defined by \eqref{eq:srg-curve} and $S_\mathrm{enc}$ to be the region enclosed by $S$.
By Theorem 4.5 of \cite{ryu2019scaled} and the arc property of the averaged operators class, we have
\[
\mathcal{G}(\cN_{\theta_1}\cN_{\theta_2})=
\Avc(\theta_1)\Avc(\theta_2).
\]
Therefore, it remains to show $\Avc(\theta_1)\Avc(\theta_2)=S_\mathrm{enc}$ on the complex plane.
The proof is completed in 3 steps.
In Step 1, we show that $S$ is the boundary enclosing $\C(\theta_1)\C(\theta_2)$ with geometric arguments.
In Step 2, we show that $S$ furthermore encloses $\Avc(\theta_1)\Avc(\theta_2)$, i.e. we show 
$\Avc(\theta_1)\Avc(\theta_2)\subseteq S_\mathrm{enc}$.
In Step 3, we show $\Avc(\theta_1)\Avc(\theta_2)= S_\mathrm{enc}$ with a topological argument.



\textbf{Step 1.}
The curve $\C(\theta_1)$ is defined by $f_1(z)=0$ with
\[f_1(z)=(x-(1-\theta_1))^2+y^2-\theta_1^2,\]
where $z=x+yi$. Let
\[
z_2(t)=1-\theta_2+\theta_2\cos(t)+\theta_2\sin(t)i
\]
be a parameterization of $\C(\theta_2)$.

Scaling and rotating $\C(\theta_1)$ by $z_2(t)\in \C(\theta_2)$ yields the curve defined by
\begin{align*}
0=&f_1(z/z_2(t))\\
=&\left(\frac{x(1-\theta_2+\theta_2\cos(t))+y\theta_2\sin(t)}{|z_2(t)|^2}-(1-\theta_1)\right)^2\\
&\hspace{1.5in}+\left(\frac{y(1-\theta_2+\theta_2\cos(t))-x\theta_2\sin(t)}{|z_2(t)|^2}\right)^2-\theta_1^2.
\end{align*}
Multiply both sides of the equation by $|z_2(t)|^4$ and simplify to get
\begin{align*}
&(x-(1-\theta_2+\theta_2\cos(t))(1-\theta_1))^2+(y+\theta_2\sin(t)(1-\theta_1))^2\\
&\hspace{1.5in}-(2\theta_2^2-2\theta_2+1+2\theta_2(1-\theta_2)\cos(t))\theta_1^2=0.
\end{align*}
Apply the envelope formula \eqref{eq:envelope_formula} to eliminate $t$ and obtain the envelope
\begin{align}
 f_2(z)&=(x^2+y^2-2x(1-\theta_1)(1-\theta_2)+(1-2\theta_1)(1-2\theta_2))^2-4\theta_1^2\theta_2^2(x^2+y^2)
\nonumber
\\
&=0.
\label{eq:f2}
\end{align}
Using polar coordinates with $r=\sqrt{x^2+y^2}\geq0$ and $x=r\cos(\varphi)$, we can factor \eqref{eq:f2} as
\begin{align*}
0=
&\big(r^2-2r\cos(\varphi)(1-\theta_1)(1-\theta_2)+(1-2\theta_1)(1-2\theta_2)- 2\theta_1\theta_2r\big)\\
&\cdot\big(r^2-2r\cos(\varphi)(1-\theta_1)(1-\theta_2)+(1-2\theta_1)(1-2\theta_2)+ 2\theta_1\theta_2r\big),
\end{align*}
where $r\ge 0$ and $\varphi\in\reals$. By considering the substitution $r\mapsto -r$ and $\varphi\mapsto \varphi+\pi$, we can combine the two factors into one to get \eqref{eq:srg-curve}:
\[
r^2-2r\cos(\varphi)(1-\theta_1)(1-\theta_2)+(1-2\theta_1)(1-2\theta_2)-2\theta_1\theta_2r=0,
\]
where $r\in \reals$ and $\varphi\in\reals$.
To clarify, the combined equation allows negative $r$.
The envelope contains contains the boundary of $\C(\theta_1)\C(\theta_2)$.

The curve defined by \eqref{eq:srg-curve} is an instance of the \emph{Cartesian oval}, which contains at most two closed curves one enclosing the other \cite{curve_book_1972}.
The following figure illustrates the envelope in solid lines and $\C(\theta_1)\C(\theta_2)$ as the shaded region.
The outer curve $S$ encloses $\C(\theta_1)\C(\theta_2)$, i.e., $\C(\theta_1)\C(\theta_2)\subseteq S_\mathrm{enc}$.

\begin{center}
		\begin{tabular}{c}
			\raisebox{-.5\height}{
				\begin{tikzpicture}[scale=2.5]
				\def\a{1/4};
				\def\b{3/4};
				\def\ti{-20};
				\def\to{-20};
				\def\tg{60};
				

				\draw[fill=medgrey,samples=200,smooth] plot[domain=-180:180]({(cos(\x)*(1-\a)*(1-\b)+\a*\b+sqrt(pow(cos(\x)*(1-\a)*(1-\b)+\a*\b,2)-(1-2*\a)*(1-2*\b)))*cos(\x)},{(cos(\x)*(1-\a)*(1-\b)+\a*\b+sqrt(pow(cos(\x)*(1-\a)*(1-\b)+\a*\b,2)-(1-2*\a)*(1-2*\b)))*sin(\x)});

				\draw[dashed,samples=200,smooth] plot[domain=-180:180]({(1-\b+\b*cos(\tg))*(1-\a+\a*cos(\x))-\a*\b*sin(\tg)*sin(\x)},{(1-\b+\b*cos(\tg))*\a*sin(\x)+(1-\a+\a*cos(\x))*\b*sin(\tg)});

				\draw[fill=white,samples=200,smooth] plot[domain=-180:180]({(cos(\x)*(1-\a)*(1-\b)-\a*\b+sqrt(pow(cos(\x)*(1-\a)*(1-\b)-\a*\b,2)-(1-2*\a)*(1-2*\b)))*cos(\x)},{(cos(\x)*(1-\a)*(1-\b)-\a*\b+sqrt(pow(cos(\x)*(1-\a)*(1-\b)-\a*\b,2)-(1-2*\a)*(1-2*\b)))*sin(\x)});

				\draw (1.5,-0.5) node [] {outer curve $S$};
				\draw [->] (1.4,-0.45) -- ({(cos(\to)*(1-\a)*(1-\b)+\a*\b+sqrt(pow(cos(\to)*(1-\a)*(1-\b)+\a*\b,2)-(1-2*\a)*(1-2*\b)))*cos(\to)},{(cos(\to)*(1-\a)*(1-\b)+\a*\b+sqrt(pow(cos(\to)*(1-\a)*(1-\b)+\a*\b,2)-(1-2*\a)*(1-2*\b)))*sin(\to)});
				
				\draw (1.05,-0.15) node [right] {inner curve};
				\draw [->] (1.05,-0.15)  -- ({(cos(\ti)*(1-\a)*(1-\b)-\a*\b+sqrt(pow(cos(\ti)*(1-\a)*(1-\b)-\a*\b,2)-(1-2*\a)*(1-2*\b)))*cos(\ti)},{(cos(\ti)*(1-\a)*(1-\b)-\a*\b+sqrt(pow(cos(\ti)*(1-\a)*(1-\b)-\a*\b,2)-(1-2*\a)*(1-2*\b)))*sin(\ti)});

				\draw (1.25,0.58) node [] {$z_2\C(\theta_1)$};
				\def\te{-65};
				\draw [->] (1.15,0.5)  -- ({(1-\b+\b*cos(\tg))*(1-\a+\a*cos(\te))-\a*\b*sin(\tg)*sin(\te)},{(1-\b+\b*cos(\tg))*\a*sin(\te)+(1-\a+\a*cos(\te))*\b*sin(\tg)});
				
				\draw [<->] (-.65,0) -- (1.2,0);
				\draw [<->] (0,-.85) -- (0,.85);
				
				\draw (1,0)  node [above right] {1};
				\filldraw (1,0) circle ({0.5*1.5/2.5pt});
				\end{tikzpicture}}
		\end{tabular}
	\end{center}


	
	



\textbf{Step 2.}
We now show that $S$ encloses not only $\C(\theta_1)\C(\theta_2)$ but also $\Avc(\theta_1)\Avc(\theta_2)$.
Note $\Avc(\theta_1)\Avc(\theta_2)$ is compact as it is the image of a compact set under a continuous map.
On the other hand, $ \Avc(\theta_1)^\circ \Avc(\theta_2)$ and $ \Avc(\theta_1) \Avc(\theta_2)^\circ$ are open as they are unions of open sets.
Since
\[
\Avc(\theta_1)^\circ \Avc(\theta_2) \cup \Avc(\theta_1) \Avc(\theta_2)^\circ
\]
is open, we have
\begin{align*}
\partial\left(\Avc(\theta_1)\Avc(\theta_2)\right)
&\subseteq
\Avc(\theta_1)\Avc(\theta_2)\backslash\left(\Avc(\theta_1)^\circ \Avc(\theta_2) \cup \Avc(\theta_1) \Avc(\theta_2)^\circ\right)\\
&\subseteq \C(\theta_1)\C(\theta_2).
\end{align*}
Since $S$ encloses $\C(\theta_1)\C(\theta_2)$ which contains the boundary of the compact set $\Avc(\theta_1)\Avc(\theta_2)$, $S$ encloses $\Avc(\theta_1)\Avc(\theta_2)$.




\textbf{Step 3.}
We have shown $\Avc(\theta_1)\Avc(\theta_2)\subseteq S_\mathrm{enc}$, and it remains to show $\Avc(\theta_1)\Avc(\theta_2)=S_\mathrm{enc}$.
The question is whether $\Avc(\theta_1)\Avc(\theta_2)$ is simply connected, i.e., whether it contains any ``holes''.
As the previous figure illustrates, $\C(\theta_1)\C(\theta_2)$ contain holes.
We show $\Avc(\theta_1)\Avc(\theta_2)$ does not.

Define the map
\begin{gather*}
\Pi\colon \Avc(\theta_1)\times \Avc(\theta_2)\rightarrow \Avc(\theta_1) \Avc(\theta_2)\\
(z_1,z_2)\mapsto z_1z_2.
\end{gather*}
To clarify, $\Avc(\theta_1)\times \Avc(\theta_2)$ denotes the product set while $\Avc(\theta_1) \Avc(\theta_2)$ denotes the Minkowski product.
We have shown that there is a parameterized closed curve 
\[
\{(z_1(t),z_2(t))\}_{t\in[0,1]}\subseteq\Avc(\theta_1)\times \Avc(\theta_2).
\]
such that $\{\Pi(z_1(t),z_2(t))\}_{t\in[0,1]}=S$.
Assume for contradiction that $z\in S_\mathrm{enc}$ but $z\notin \Avc(\theta_1)\Avc(\theta_2)$. (In other words, we assume for contradiction that $z$ is strictly within the hole of the domain.) Since $\Avc(\theta_1)\times \Avc(\theta_2)$ is simply connected, we can continuously contract $\{(z_1(t),z_2(t))\}_{t\in[0,1]}$ to a point in $\Avc(\theta_1)\times \Avc(\theta_2)$, and the curve under the map $\Pi$ continuously contracts to a point in $\Avc(\theta_1) \Avc(\theta_2)$.
However, this is not possible as $\{\Pi(z_1(t),z_2(t))\}_{t\in[0,1]}$ has a nonzero winding number around $z$ and $z\notin \Avc(\theta_1) \Avc(\theta_2)$.
We have a contradiction and we conclude $z\in \Avc(\theta_1) \Avc(\theta_2)$.
\qed


\subsection{Proof of Corollary~\ref{cor:ogura}}
We can visually observe from Figure~\ref{fig:srg-exmaples}
that $\C(\theta)$, the dashed circle, enclose $\mathcal{G}(\cN_{\theta_1}\cN_{\theta_2})$.
We can also observe that the geometric objects have matching curvature at point $1$, and therefore we cannot further reduce the size of the dashed circle while enclosing $\mathcal{G}(\cN_{\theta_1}\cN_{\theta_2})$. We now make this argument formal with Fact~\ref{fact:srg} and the following geometric arguments.

Remember, $f_2(x,y)=0$ defines the boundary $\partial(\Avc(\theta_1)\Avc(\theta_2))$.
Define
\begin{align*}
g(\varphi)&=
f_2(\theta\cos(\varphi)+(1-\theta),\theta\sin(\varphi))\\
&=
\frac{16  \theta_1^2\theta_2^2 (1-\theta_1)^2(1-\theta_2)^2 (\theta_1+\theta_2-2 \theta_1 \theta_2)^2 \sin ^4(\varphi/2)}{(1-\theta_1 \theta_2)^4},
\end{align*}
i.e., $g(\varphi)$ is $f_2$ evaluated on the curve $\C(\theta)$.
We can see that $g(\varphi)>0$ for all $\varphi\ne 0$ and $g(0)=0$.
This implies $\C(\theta)$ and $\partial(\Avc(\theta_1)\Avc(\theta_2))$ intersect at only one point and therefore do not cross.
The point $(1-\varepsilon,0)$ is in $\Avc(\theta_1)\Avc(\theta_2)$ and enclosed by $\C(\theta)$ for small enough $\varepsilon>0$.
Since
\[
f_2(1-\varepsilon,0)=-8\theta_1\theta_2\underbrace{(\theta_1+\theta_2-2\theta_1\theta_2)}_{>0}\varepsilon+\mathcal{O}(\varepsilon^2)
\]
for $\varepsilon\rightarrow 0$, it is $\C(\theta)$ that encloses $\partial(\Avc(\theta_1)\Avc(\theta_2))$.
Finally, we conclude $\Avc(\theta)$ contains $\Avc(\theta_1)\Avc(\theta_2)$.

Consider $r(\varphi)$ defined by \eqref{eq:srg-curve}.
Through implicit differentiation, we get
\[
\frac{dr}{d\varphi}\big|_{\varphi=0}=0,\qquad\frac{d^2r}{d\varphi^2}\big|_{\varphi=0}=\frac{(1-\theta_1)(1-\theta_2)}{\theta_1+\theta_2-2\theta_1\theta_2}.
\]
Using \eqref{eq:curvature}, the curvature of $r(\varphi)$ at point $1$ (given by $\varphi=0$) is 
\[
\kappa(\varphi)\big|_{\varphi=0}=\frac{r(\varphi)^2+2(\frac{dr}{d\varphi})^2-r(\varphi)\frac{d^2r}{d\varphi^2}}{(r(\varphi)^2+(\frac{dr}{d\varphi})^2)^{\frac{3}{2}}}\big|_{\varphi=0}=\frac{1-\theta_1\theta_2}{\theta_1+\theta_2-2\theta_1\theta_2}=\frac{1}{\theta}.
\]
This implies any circle through 1 symmetric about the real axis containing $r(\varphi)$ must have radius at least $\theta$.
\qed

\section{Tight characterization of Davis--Yin splitting}
\label{s:dys}
Consider the monotone inclusion problem
\[
\begin{array}{ll}
\underset{x\in \hilbert}{\mbox{find}}&0\in (A+B+C)x,
\end{array}
\]
where $A$ and $B$ are maximal monotone and $C$ is $\beta$-cocoercive.
Davis and Yin \cite{davis2017three} proposed
\[
T_{\gamma}(A,B,C)=I-J_{\gamma B}+J_{\gamma A} (2J_{\gamma B}-I-\gamma C  J_{\gamma B}),
\]
which we call the \emph{Davis--Yin splitting} (DYS),
and showed that it is a fixed-point encoding for the monotone inclusion problem in the sense that
$\mathrm{Zer}(A+B+C)=J_{\gamma B}(\fix (T_{\gamma}(A,B,C)))$,
where $\mathrm{Zer}$ and $\fix$ repsectively denote the set of zeros and fixed points.
Define the class of DYS operators as
\[
\cT_{\beta,\gamma}=\left\{T_{\gamma}(A,B,C)\,\big|\,A,B\in\cM,\,C\in\cC_{\beta}\right\}.
\]
Davis and Yin showed that the DYS operators of $\cT_{\beta,\gamma}$ are $\frac{2\beta}{4\beta-\gamma}$-averaged.
\begin{fact}[Proposition~2.1 of \cite{davis2017three}]
\label{prop:davis}
Let $\gamma\in(0,2\beta)$. Then

\vspace{-0.15in}
\[
\cG(\cT_{\beta,\gamma})\subseteq\cG\left(\cN_{\frac{2\beta}{4\beta-\gamma}}\right)=
\hspace{-0.1in}
\begin{tabular}{c}
\raisebox{-.5\height}{
\begin{tikzpicture}[scale=1.8]

\def\b{1};
\def\c{1};
\def\g{1.3};
\def\gg{1};
\def\t{18};

\coordinate (O) at (0,0);
\coordinate (A) at ({cos(\t)^2},{cos(\t)*sin(\t)});
\coordinate (X) at (1,0);
\coordinate (I) at ({1-2\b*\c/(4-\g)},0);
\coordinate (A_2) at ({1-\b*\c+\b*\c*cos(2*\t)^2},{\b*\c*cos(2*\t)*sin(2*\t)});
\coordinate (O_4) at ({1-\b*\c+\b*\c*(cos(2*\t)-\g/2*cos(\t)^2)*cos(2*\t)},{\b*\c*(cos(2*\t)-\g/2*cos(\t)^2)*sin(2*\t)});
\coordinate (A_3) at ({\g*\b*\c*cos(\t)^2*(cos(2*\t))},{\g*\b*\c*cos(\t)^2*(sin(2*\t))});
\coordinate (O_3) at ({0.5*\g*\b*\c*cos(\t)^2*(cos(2*\t))},{0.5*\g*\b*\c*cos(\t)^2*(sin(2*\t))});
\coordinate (I) at ({1-2*\b*\c/(4-\g)},0);
\coordinate (II) at ({1-2*\b*\c/(4-\gg)},0);
\coordinate (B) at (-0.13,0.63);
\coordinate (T) at (0.99,0.41);
\coordinate (T') at (-0.34,0.06);

\fill [medgrey] (I) circle ({2*\b*\c/(4-\g)});

\draw [<->] (-0.6,0) -- (1.2,0);
\draw [<->] (0,-0.8) -- (0,0.8);

\filldraw (X) circle[radius={0.5*1.5/1.8pt}];
\filldraw (I) circle[radius={0.5*1.5/1.8pt}];

\draw (X) node[above right] {$1$};
\draw (0.25,0) node[below] {$\frac{2\beta-\gamma}{4\beta-\gamma}$};


\end{tikzpicture}}
\end{tabular}
\]
\end{fact}

We show that this characterization is tight in the following sense.
\begin{theorem}
\label{thm:dys-srg}
Let $\gamma\in(0,2\beta)$. Then
\vspace{-0.3in}

\[
\cG(\cT_{\beta,\gamma})=\cG\left(\cN_{\frac{2\beta}{4\beta-\gamma}}\right)=
\hspace{-0.1in}
\begin{tabular}{c}
\raisebox{-.5\height}{
\begin{tikzpicture}[scale=1.8]

\def\b{1};
\def\c{1};
\def\g{1.3};
\def\gg{1};
\def\t{18};

\coordinate (O) at (0,0);
\coordinate (A) at ({cos(\t)^2},{cos(\t)*sin(\t)});
\coordinate (X) at (1,0);
\coordinate (I) at ({1-2\b*\c/(4-\g)},0);
\coordinate (A_2) at ({1-\b*\c+\b*\c*cos(2*\t)^2},{\b*\c*cos(2*\t)*sin(2*\t)});
\coordinate (O_4) at ({1-\b*\c+\b*\c*(cos(2*\t)-\g/2*cos(\t)^2)*cos(2*\t)},{\b*\c*(cos(2*\t)-\g/2*cos(\t)^2)*sin(2*\t)});
\coordinate (A_3) at ({\g*\b*\c*cos(\t)^2*(cos(2*\t))},{\g*\b*\c*cos(\t)^2*(sin(2*\t))});
\coordinate (O_3) at ({0.5*\g*\b*\c*cos(\t)^2*(cos(2*\t))},{0.5*\g*\b*\c*cos(\t)^2*(sin(2*\t))});
\coordinate (I) at ({1-2*\b*\c/(4-\g)},0);
\coordinate (II) at ({1-2*\b*\c/(4-\gg)},0);
\coordinate (B) at (-0.13,0.63);
\coordinate (T) at (0.99,0.41);
\coordinate (T') at (-0.34,0.06);

\fill [medgrey] (I) circle ({2*\b*\c/(4-\g)});

\draw [<->] (-0.6,0) -- (1.2,0);
\draw [<->] (0,-0.8) -- (0,0.8);

\filldraw (X) circle[radius={0.5*1.5/1.8pt}];
\filldraw (I) circle[radius={0.5*1.5/1.8pt}];

\draw (X) node[above right] {$1$};
\draw (0.25,0) node[below] {$\frac{2\beta-\gamma}{4\beta-\gamma}$};


\end{tikzpicture}}
\end{tabular}
\]
\end{theorem}
\begin{corollary}
\label{cor:dys}
The averagedness parameter of Fact~\ref{prop:davis} is tight in the sense that it cannot be improved without further assumptions.
\end{corollary}

\subsection{Proof of Theorem~\ref{thm:dys-srg}}
Since $\cG(\cT_{\beta,\gamma})\subseteq \Avc\left(\frac{2\beta}{4\beta-\gamma}\right)$ by Fact~\ref{prop:davis}, we show $\Avc\left(\frac{2\beta}{4\beta-\gamma}\right)\subseteq \cG(\cT_{\beta,\gamma})$.
Define the set
\[
S_{\beta,\gamma}=\left\{1-z_2+z_1(2z_2-1-\gamma z_3z_2)\,|\,z_1,z_2\in\Avc(1/2),\,z_3\in\tfrac{1}{\beta}\Avc(1/2)\right\}.
\]
The proof is completed in three steps.
In Step 1, we show $S_{\beta,\gamma}\subseteq \cG(\cT_{\beta,\gamma})$ by appealing to results about the SRG.
In Step 2, we show $\C\left(\frac{2\beta}{4\beta-\gamma}\right)\subseteq S_{\beta,\gamma}$ with geometric arguments.
In Step 3, we strengthen the result of Step 2 to $\Avc\left(\frac{2\beta}{4\beta-\gamma}\right)\subseteq S_{\beta,\gamma}$ using a topological argument.

\textbf{Step 1.}
By Lemma 3.2, Proposition 3.3, and Theorems 4.2 and 4.3 of \cite{ryu2019scaled}, we can identify $z_1,z_2\in\Avc(1/2)$ with resolvents of maximal monotone operators on $\reals^2$ and $z_3\in\tfrac{1}{\beta}\Avc(1/2)$ with a $\beta$-cocoercive operator on $\reals^2$.
Therefore, $S_{\beta,\gamma}$ represents the SRGs of operators in $\cT_{\beta,\gamma}$, and we conclude $S_{\beta,\gamma}\subseteq \cG(\cT_{\beta,\gamma})$.



\textbf{Step 2.}
Define
\[
R_{\beta,\gamma}=\left\{1-z_2+z_1(2z_2-1-\gamma z_3z_2)\,\Big|\,z_1=z_2\in\Avc(1/2),\,z_3\in\tfrac{1}{\beta}\Avc(1/2)\right\}.
\]
Clearly $R_{\beta,\gamma}\subseteq S_{\beta,\gamma}$.
We show
\begin{align*}
\C\left(\frac{2\beta}{4\beta-\gamma}\right)\subseteq R_{\beta,\gamma}.
\end{align*}

Let $A_1= z_1=z_2=\cos(\theta)e^{i\theta} \in\C(1/2)$.
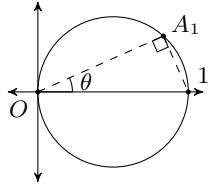
\begin{center}
\begin{tabular}{c}
\raisebox{-.5\height}{
\begin{tikzpicture}[scale=2]
\def\t{24};

\coordinate (O) at (0,0);
\coordinate (A) at ({cos(\t)*cos(\t)},{cos(\t)*sin(\t)});
\coordinate (X) at (1,0);

\draw [<->] (-0.2,0) -- (1.1,0);
\draw [<->] (0,-0.6) -- (0,0.6);

\filldraw (A) circle[radius={0.5*1.5/2pt}];
\filldraw (X) circle[radius={0.5*1.5/2pt}];
\filldraw (O) circle[radius={0.5*1.5/2pt}];

\draw (1/2,0) circle (1/2);

\draw [dashed] (O) -- (A);
\draw [dashed] (A) -- (X);

\draw (X) node[above right] {$1$};
\draw  ({cos(\t)*cos(\t)-0.01},{cos(\t)*sin(\t)-0.05}) node[above right] {$A_1$};
\draw (0.22,0.07) node[right] {$\theta$};
\draw (O) node[below left] {$O$};

\draw [] (0,0) -- (0.23,0) arc (0:\t:0.23);

\tkzMarkRightAngle[size=.08](O,A,X);

\end{tikzpicture}}
\end{tabular}
\end{center}
With direct calculations, we have 
\[
A_2= 2z_1z_2-z_1-z_2+1=\cos(2\theta)e^{2\theta i}\in \C(1/2)
\]
and 
\[
A_3= \frac{\gamma}{\beta} z_1z_2=\frac{\gamma}{\beta} \cos^2(\theta)e^{2\theta i}.
\]
Define $O_1=\frac{A_3}{2}$.
Figure~\ref{fig:step2_illustration} illustrates the following construction.
Define point $P=\frac{2\beta-\gamma}{4\beta-\gamma}$ as the center of the circle $\C\left(\frac{2\beta}{4\beta-\gamma}\right)$.
Let $O_2=A_2-\frac{A_3}{2}=A_2-O_1$ be the center of the disk $A_2-A_3\Avc(1/2)$.
Let $B$ be the point farthest from $P$ in the disk $A_2-A_3\Avc(1/2)$.
($B\in R_{\beta,\gamma}$ since $A_2-A_3\Avc(1/2)\subseteq R_{\beta,\gamma}$.)
Then $P$, $O_2$, and $B$ are collinear. We have 
\[
\overline{O_2B}=\frac{\gamma}{2\beta}\cos^2(\theta).
\]
Using the cosine rule, we have
\begin{align*}
\overline{PO_2}^2
&=
\overline{O_2O}^2+\overline{OP}^2-2 
\overline{O_2O}\cdot\overline{OP}\cos(2\theta)\\
&=
{\scriptstyle \left(\cos(2\theta)-\frac{\gamma}{2\beta}\cos^2(\theta)\right)^2
+\left(\frac{2\beta-\gamma}{4\beta-\gamma}\right)^2-
2 
\left(\cos(2\theta)-\frac{\gamma}{2\beta}\cos^2(\theta)\right)
\left(\frac{2\beta-\gamma}{4\beta-\gamma}\right)
\cos(2\theta)}\\
&=
\left(\frac{\gamma}{2\beta}\cos^2(\theta)-\frac{2\beta}{4\beta-\gamma}\right)^2.
\end{align*}
Since $P$, $O_2$, and $B$ are collinear, we have
\[
\overline{PB}=\overline{PO_2}+\overline{O_2B}=\frac{2\beta}{4\beta-\gamma}.
\]
Therefore, $B\in \C\left(\frac{2\beta}{4\beta-\gamma}\right)$.
	\newpage
	
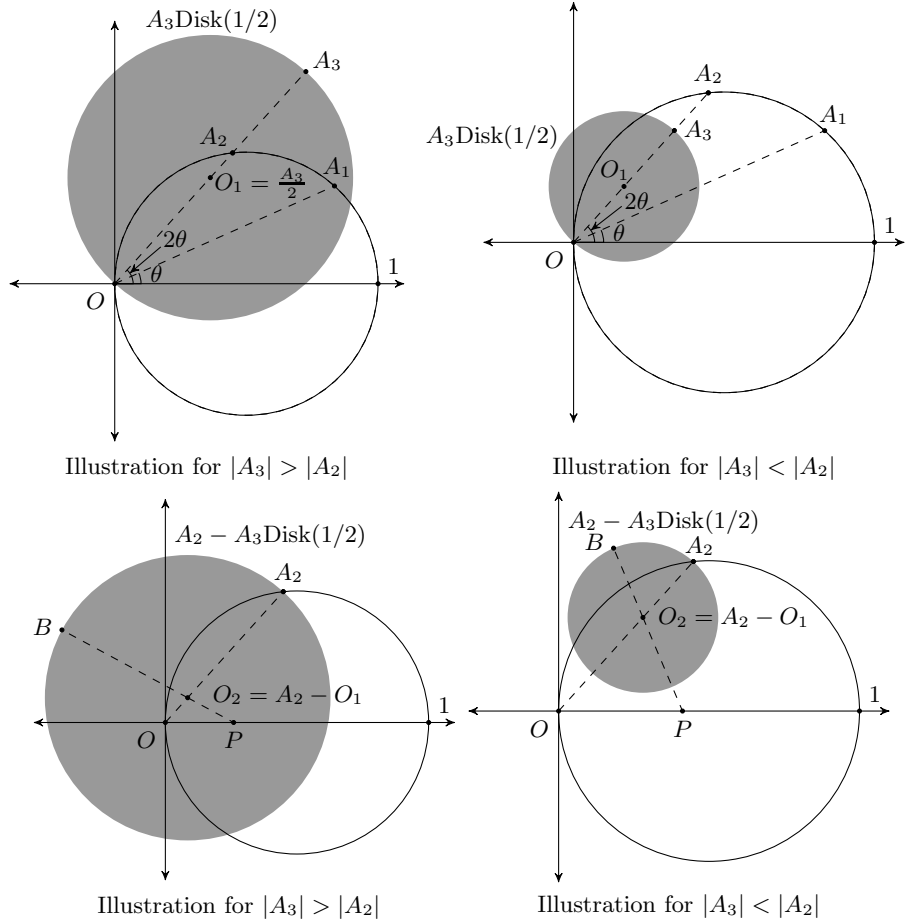
\begin{figure}
\centering
			\begin{tikzpicture}[scale=3.5]
			
			\def\b{1};
			\def\c{1};
			\def\g{1.3};
			\def\t{24};

			\coordinate (O) at (0,0);
			\coordinate (A) at ({cos(\t)^2},{cos(\t)*sin(\t)});
			\coordinate (X) at (1,0);
			\coordinate (I) at ({1-2\b*\c/(4-\g)},0);
			\coordinate (A_2) at ({1-\b*\c+\b*\c*cos(2*\t)^2},{\b*\c*cos(2*\t)*sin(2*\t)});
			\coordinate (O_4) at ({1-\b*\c+\b*\c*(cos(2*\t)-\g/2*cos(\t)^2)*cos(2*\t)},{\b*\c*(cos(2*\t)-\g/2*cos(\t)^2)*sin(2*\t)});
			\coordinate (A_3) at ({\g*\b*\c*cos(\t)^2*(cos(2*\t))},{\g*\b*\c*cos(\t)^2*(sin(2*\t))});
			\coordinate (O_3) at ({0.5*\g*\b*\c*cos(\t)^2*(cos(2*\t))},{0.5*\g*\b*\c*cos(\t)^2*(sin(2*\t))});
			\coordinate (I) at ({1-2*\b*\c/(4-\g)},0);
			\coordinate (B) at (-0.13,0.63);
			\coordinate (T) at (0.99,0.41);
			\coordinate (T') at (-0.34,0.06);

			\fill [medgrey] (O_3) circle ({0.5*\g*\b*\c*cos(\t)^2});

			\draw (1/2,0) circle (1/2);
			\draw [dashed] ({1-\b*\c/2},0) circle ({\b*\c/2});
			
			\draw [<->] (-0.4,0) -- (1.1,0);
			\draw [<->] (0,-0.6) -- (0,1);
			
			\filldraw (A) circle[radius={0.5*1.5/3.5pt}];
			\filldraw (X) circle[radius={0.5*1.5/3.5pt}];
			\filldraw (A_2) circle[radius={0.5*1.5/3.5pt}];
			\filldraw (A_3) circle[radius={0.5*1.5/3.5pt}];
			\filldraw (O) circle[radius={0.5*1.5/3.5pt}];
			\filldraw (O_3) circle[radius={0.5*1.5/3.5pt}];
			\filldraw (A_2) circle[radius={0.5*1.5/3.5pt}];
			
			\draw [dashed] (O) -- (A);
			\draw [dashed] (O) -- (A_3);

			\draw (X) node[above right] {$1$};
			\draw ({cos(\t)^2+0.01},{cos(\t)*sin(\t)}) node[above] {$A_1$};
			\draw ({\g*\b*\c*cos(\t)^2*(cos(2*\t))-0.01},{\g*\b*\c*cos(\t)^2*(sin(2*\t))-0.03}) node[above right] {$A_3$};

			\draw ({1-\b*\c+\b*\c*cos(2*\t)^2+0.02},{\b*\c*cos(2*\t)*sin(2*\t)-0.01}) node[above left] {$A_2$};
			\draw (O) node[below left] {$O$};
			\draw ({0.5*\g*\b*\c*cos(\t)^2*(cos(2*\t))-0.02},{0.5*\g*\b*\c*cos(\t)^2*(sin(2*\t))-0.02}) node[right] {$O_1=\frac{A_3}{2}$};

			\draw (0.37,1.0) node {$A_3\Avc(1/2)$};

			\draw [] (0,0) -- (0.1,0) arc (0:\t:0.1);
			\draw (0.1,0.038) node[right] {$\theta$};

			\draw [] (0,0) -- (0.07,0) arc (0:{2*\t}:0.07);
			\draw (0.23,0.17) node[] {$2\theta$};
			\draw [->] (0.175,0.14) -- ({0.07*cos(35)},{0.07*sin(35)});
			
			\draw (0.35,-0.7) node {Illustration for $|A_3|>|A_2|$};
			
			\end{tikzpicture}
			\begin{tikzpicture}[scale=4]
			
			\def\b{1};
			\def\c{1};
			\def\g{0.6};
			\def\t{24};
			
			\coordinate (O) at (0,0);
			\coordinate (A) at ({cos(\t)^2},{cos(\t)*sin(\t)});
			\coordinate (X) at (1,0);
			\coordinate (I) at ({1-2\b*\c/(4-\g)},0);
			\coordinate (A_2) at ({1-\b*\c+\b*\c*cos(2*\t)^2},{\b*\c*cos(2*\t)*sin(2*\t)});
			\coordinate (O_4) at ({1-\b*\c+\b*\c*(cos(2*\t)-\g/2*cos(\t)^2)*cos(2*\t)},{\b*\c*(cos(2*\t)-\g/2*cos(\t)^2)*sin(2*\t)});
			\coordinate (A_3) at ({\g*\b*\c*cos(\t)^2*(cos(2*\t))},{\g*\b*\c*cos(\t)^2*(sin(2*\t))});
			\coordinate (O_3) at ({0.5*\g*\b*\c*cos(\t)^2*(cos(2*\t))},{0.5*\g*\b*\c*cos(\t)^2*(sin(2*\t))});
			\coordinate (I) at ({1-2*\b*\c/(4-\g)},0);
			\coordinate (B) at (-0.13,0.63);
			\coordinate (T) at (0.99,0.41);
			\coordinate (T') at (-0.34,0.06);

			\fill [medgrey] (O_3) circle ({0.5*\g*\b*\c*cos(\t)^2});

			\draw (1/2,0) circle (1/2);
			\draw [dashed] ({1-\b*\c/2},0) circle ({\b*\c/2});

			\draw [<->] (-0.3,0) -- (1.1,0);
			\draw [<->] (0,-0.6) -- (0,0.8);
			
			\filldraw (A) circle[radius={0.5*1.5/4pt}];
			\filldraw (O) circle[radius={0.5*1.5/4pt}];
			\filldraw (X) circle[radius={0.5*1.5/4pt}];
			\filldraw (A_2) circle[radius={0.5*1.5/4pt}];
			\filldraw (A_3) circle[radius={0.5*1.5/4pt}];
			\filldraw (O_3) circle[radius={0.5*1.5/4pt}];
			\filldraw (A_2) circle[radius={0.5*1.5/4pt}];

			\draw [dashed] (O) -- (A);
			\draw [dashed] (O) -- (A_2);

			\draw (X) node[above right] {$1$};
			\draw ({cos(\t)^2+0.03},{cos(\t)*sin(\t)-0.01}) node[above] {$A_1$};
			\draw (A_3) node[right] {$A_3$};

			\draw (A_2) node[above] {$A_2$};
			\draw (O) node[below left] {$O$};
			\draw ({0.5*\g*\b*\c*cos(\t)^2*(cos(2*\t))-0.03},{0.5*\g*\b*\c*cos(\t)^2*(sin(2*\t))-0.01}) node[above] {$O_1$};

			\draw (-0.27,0.35) node {$A_3\Avc(1/2)$};
			
			\draw [] (0,0) -- (0.1,0) arc (0:\t:0.1);
			\draw (0.1,0.033) node[right] {$\theta$};
			
			\draw [] (0,0) -- (0.07,0) arc (0:{2*\t}:0.07);
			\draw (0.21,0.08) node[above] {$2\theta$};
			\draw [->] (0.16,0.12) -- ({0.07*cos(35)},{0.07*sin(35)});
			
			\draw (0.4,-0.75) node {Illustration for $|A_3|<|A_2|$};
			\end{tikzpicture}\\
			\begin{tikzpicture}[scale=3.5]
			
			\def\b{1};
			\def\c{1};
			\def\g{1.3};
			\def\t{24};
			
			\def\j{(2-\g)/(4-\g)};
			\def\l{2/((4-\g)*sqrt((cos(2*\t)-\g/2*cos(\t)^2)^2+(\j)^2-2*\j*(cos(2*\t)-\g/2*cos(\t)^2)*cos(2*\t)))};

			\coordinate (O) at (0,0);
			\coordinate (A) at ({cos(\t)^2},{cos(\t)*sin(\t)});
			\coordinate (X) at (1,0);
			\coordinate (I) at ({1-2\b*\c/(4-\g)},0);
			\coordinate (A_2) at ({1-\b*\c+\b*\c*cos(2*\t)^2},{\b*\c*cos(2*\t)*sin(2*\t)});
			\coordinate (O_4) at ({1-\b*\c+\b*\c*(cos(2*\t)-\g/2*cos(\t)^2)*cos(2*\t)},{\b*\c*(cos(2*\t)-\g/2*cos(\t)^2)*sin(2*\t)});
			\coordinate (A_3) at ({\g*\b*\c*cos(\t)^2*(cos(2*\t))},{\g*\b*\c*cos(\t)^2*(sin(2*\t))});
			\coordinate (O_3) at ({0.5*\g*\b*\c*cos(\t)^2*(cos(2*\t))},{0.5*\g*\b*\c*cos(\t)^2*(sin(2*\t))});
			\coordinate (I) at ({1-2*\b*\c/(4-\g)},0);
			\coordinate (B) at ({\j+((cos(2*\t)-\g/2*cos(\t)^2)*cos(2*\t)-\j)*\l},{(cos(2*\t)-\g/2*cos(\t)^2)*sin(2*\t)*\l});

			\fill [medgrey] (O_4) circle ({0.5*\g*\b*\c*cos(\t)^2});
			\draw (1/2,0) circle (1/2);

			\draw [<->] (-0.5,0) -- (1.1,0);
			\draw [<->] (0,-0.6) -- (0,0.85);

			\filldraw (X) circle[radius={0.5*1.5/3.5pt}];
			\filldraw (A_2) circle[radius={0.5*1.5/3.5pt}];
			\filldraw (O_4) circle[radius={0.5*1.5/3.5pt}];
			\filldraw (A_2) circle[radius={0.5*1.5/3.5pt}];
			\filldraw (I) circle[radius={0.5*1.5/3.5pt}];
			\filldraw (B) circle[radius={0.5*1.5/3.5pt}];
			\filldraw (O) circle[radius={0.5*1.5/3.5pt}];

			\draw [dashed] (O) -- (A_2);
			\draw [dashed] (O_4) -- (I);
			\draw [dashed] (O_4) -- (B);
			
			\draw (X) node[above right] {$1$};

			\draw ({1-\b*\c+\b*\c*cos(2*\t)^2+0.02},{\b*\c*cos(2*\t)*sin(2*\t)}) node[above] {$A_2$};
			\draw (O) node[below left] {$O$};
			\draw ({1-\b*\c+\b*\c*(cos(2*\t)-\g/2*cos(\t)^2)*cos(2*\t)+0.06},{\b*\c*(cos(2*\t)-\g/2*cos(\t)^2)*sin(2*\t)+0.01}) node[right] {$O_2=A_2-O_1$};
			\draw (I) node[below] {$P$};
			\draw ({\j+((cos(2*\t)-\g/2*cos(\t)^2)*cos(2*\t)-\j)*\l},{(cos(2*\t)-\g/2*cos(\t)^2)*sin(2*\t)*\l+0.01}) node[left] {$B$};
			\draw (0.4,0.71) node {$A_2-A_3\Avc(1/2)$};

			
			\draw (0.28,-0.7) node {Illustration for $|A_3|>|A_2|$};
			\end{tikzpicture}
			\begin{tikzpicture}[scale=4]
			
			\def\b{1};
			\def\c{1};
			\def\g{0.6};
			\def\t{24};
			
			\def\j{(2-\g)/(4-\g)};
			\def\l{2/((4-\g)*sqrt((cos(2*\t)-\g/2*cos(\t)^2)^2+(\j)^2-2*\j*(cos(2*\t)-\g/2*cos(\t)^2)*cos(2*\t)))};
			
			\coordinate (O) at (0,0);
			\coordinate (A) at ({cos(\t)^2},{cos(\t)*sin(\t)});
			\coordinate (X) at (1,0);
			\coordinate (I) at ({1-2\b*\c/(4-\g)},0);
			\coordinate (A_2) at ({1-\b*\c+\b*\c*cos(2*\t)^2},{\b*\c*cos(2*\t)*sin(2*\t)});
			\coordinate (O_4) at ({1-\b*\c+\b*\c*(cos(2*\t)-\g/2*cos(\t)^2)*cos(2*\t)},{\b*\c*(cos(2*\t)-\g/2*cos(\t)^2)*sin(2*\t)});
			\coordinate (A_3) at ({\g*\b*\c*cos(\t)^2*(cos(2*\t))},{\g*\b*\c*cos(\t)^2*(sin(2*\t))});
			\coordinate (O_3) at ({0.5*\g*\b*\c*cos(\t)^2*(cos(2*\t))},{0.5*\g*\b*\c*cos(\t)^2*(sin(2*\t))});
			\coordinate (I) at ({1-2*\b*\c/(4-\g)},0);
			\coordinate (B) at ({\j+((cos(2*\t)-\g/2*cos(\t)^2)*cos(2*\t)-\j)*\l},{(cos(2*\t)-\g/2*cos(\t)^2)*sin(2*\t)*\l});

			\fill [medgrey] (O_4) circle ({0.5*\g*\b*\c*cos(\t)^2});
			\draw (1/2,0) circle (1/2);

			\draw [<->] (-0.3,0) -- (1.1,0);
			\draw [<->] (0,-0.57) -- (0,0.73);

			\filldraw (X) circle[radius={0.5*1.5/4pt}];
			\filldraw (A_2) circle[radius={0.5*1.5/4pt}];
			\filldraw (O_4) circle[radius={0.5*1.5/4pt}];
			\filldraw (A_2) circle[radius={0.5*1.5/4pt}];
			\filldraw (I) circle[radius={0.5*1.5/4pt}];
			\filldraw (B) circle[radius={0.5*1.5/4pt}];
			\filldraw (O) circle[radius={0.5*1.5/4pt}];

			\draw [dashed] (O) -- (A_2);
			\draw [dashed] (O_4) -- (I);
			\draw [dashed] (O_4) -- (B);

			\draw (X) node[above right] {$1$};

			\draw ({1-\b*\c+\b*\c*cos(2*\t)^2+0.02},{\b*\c*cos(2*\t)*sin(2*\t)-0.01}) node[above] {$A_2$};
			\draw (O) node[below left] {$O$};
			\draw ({1-\b*\c+\b*\c*(cos(2*\t)-\g/2*cos(\t)^2)*cos(2*\t)+0.02},{\b*\c*(cos(2*\t)-\g/2*cos(\t)^2)*sin(2*\t)}) node[right] {$O_2=A_2-O_1$};
			\draw (I) node[below] {$P$};
			\draw ({\j+((cos(2*\t)-\g/2*cos(\t)^2)*cos(2*\t)-\j)*\l},{(cos(2*\t)-\g/2*cos(\t)^2)*sin(2*\t)*\l+0.02}) node[left] {$B$};
			\draw (0.35,0.63) node {$A_2-A_3\Avc(1/2)$};
			
			
			\draw (0.4,-0.65) node {Illustration for $|A_3|<|A_2|$};
			\end{tikzpicture}
	\caption{The first geometric construction of Step 2.
	The shaded region is a subset of $R_{\beta,\gamma}$ by construction.
	}
	\label{fig:step2_illustration}
\end{figure}

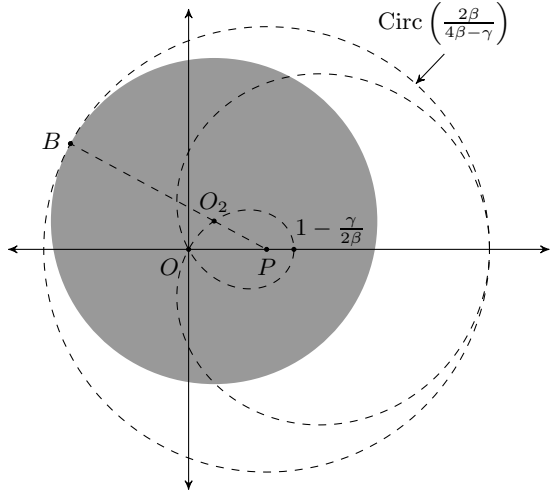
\begin{figure}
\centering
			\begin{tikzpicture}[scale=4]
			
			\def\b{1};
			\def\c{1};
			\def\g{1.3};
			\def\t{24};

			\def\j{(2-\g)/(4-\g)};
			\def\l{2/((4-\g)*sqrt((cos(2*\t)-\g/2*cos(\t)^2)^2+(\j)^2-2*\j*(cos(2*\t)-\g/2*cos(\t)^2)*cos(2*\t)))};
			
			\coordinate (O) at (0,0);
			\coordinate (A) at ({cos(\t)^2},{cos(\t)*sin(\t)});
			\coordinate (X) at (1,0);
			\coordinate (I) at ({1-2\b*\c/(4-\g)},0);
			\coordinate (A_2) at ({1-\b*\c+\b*\c*cos(2*\t)^2},{\b*\c*cos(2*\t)*sin(2*\t)});
			\coordinate (O_4) at ({1-\b*\c+\b*\c*(cos(2*\t)-\g/2*cos(\t)^2)*cos(2*\t)},{\b*\c*(cos(2*\t)-\g/2*cos(\t)^2)*sin(2*\t)});
			\coordinate (A_3) at ({\g*\b*\c*cos(\t)^2*(cos(2*\t))},{\g*\b*\c*cos(\t)^2*(sin(2*\t))});
			\coordinate (O_3) at ({0.5*\g*\b*\c*cos(\t)^2*(cos(2*\t))},{0.5*\g*\b*\c*cos(\t)^2*(sin(2*\t))});
			\coordinate (I) at ({1-2*\b*\c/(4-\g)},0);
			\coordinate (B) at ({\j+((cos(2*\t)-\g/2*cos(\t)^2)*cos(2*\t)-\j)*\l},{(cos(2*\t)-\g/2*cos(\t)^2)*sin(2*\t)*\l});

			\fill [medgrey] (O_4) circle ({0.5*\g*\b*\c*cos(\t)^2});
			\draw [dashed] (I) circle ({2*\b*\c/(4-\g)});
			\draw [dashed,samples=200,smooth] plot[domain=-180:180] (\x:{cos(\x)-0.5*\g*cos(0.5*\x)*cos(0.5*\x)});
			
			\draw [<->] (-0.6,0) -- (1.2,0);
			\draw [<->] (0,-0.8) -- (0,0.8);
			
			\filldraw (O) circle[radius={0.5*1.5/4pt}];
			\filldraw (O_4) circle[radius={0.5*1.5/4pt}];
			\filldraw (I) circle[radius={0.5*1.5/4pt}];
			\filldraw (B) circle[radius={0.5*1.5/4pt}];
			\filldraw ({1-\g/2},0) circle[radius={0.5*1.5/4pt}];

			\draw [dashed] (O_4) -- (I);
			\draw [dashed] (O_4) -- (B);

			\draw (O) node[below left] {$O$};
			\draw (O_4) node[above] {$O_2$};
            \draw (I) node[below] {$P$};
			\draw (0.47,-0.02) node[above] {$1-\frac{\gamma}{2\beta}$};
			\draw ({\j+((cos(2*\t)-\g/2*cos(\t)^2)*cos(2*\t)-\j)*\l},{(cos(2*\t)-\g/2*cos(\t)^2)*sin(2*\t)*\l+0.01}) node[left] {$B$};
			
			\draw (0.85,0.65) node[above] {$\C\left(\frac{2\beta}{4\beta-\gamma}\right)$};
			\draw [->] (0.85,0.65) -- ({1-2*\b*\c/(4-\g)+2*\b*\c/(4-\g)*cos(48)},{2*\b*\c/(4-\g)*sin(48)});
			\end{tikzpicture}
				\caption{The second geometric construction of Step 2.}
			\label{fig:Step2_step2}
\end{figure}

Figure~\ref{fig:Step2_step2} illustrates the following construction.
The trajectory of $O_2=(\cos(2\theta)-\frac{\gamma}{2\beta} \cos^2(\theta))e^{2\theta i}$ as a function of $\theta$ is a closed curve within $\C\left(\frac{2\beta}{4\beta-\gamma}\right)$.
Since $0<1-\frac{2\beta}{4\beta-\gamma}<1-\frac{\gamma}{2\beta}$, the curve strictly encloses $P$.
As $\theta$ traverses $[-\pi/2,\pi/2)$, $O_2$ traverses the inner curve and $B$ traverses all of $\C\left(\frac{2\beta}{4\beta-\gamma}\right)$.
Therefore, we conclude	$\C\left(\frac{2\beta}{4\beta-\gamma}\right)\subseteq S_{\beta,\gamma}$.

\textbf{Step 3.}
Define the map 
\begin{gather*}
\Pi\colon \Avc(1/2)\times \Avc(1/2) \times\tfrac{1}{\beta}\Avc(1/2)\rightarrow S_{\beta,\gamma}\\
(z_1,z_2,z_3)\mapsto 1-z_2+z_1(2z_2-1-\gamma z_3z_2).
\end{gather*}
Consider any $z$ strictly enclosed within $\C\left(\frac{2\beta}{4\beta-\gamma}\right)$, and assume for contradiction that $z\notin S_{\beta,\gamma}$.
We have shown that there is a closed curve 
\[
\{\eta(t)\}_{t\in[0,1]}\subseteq \Avc(1/2)\times \Avc(1/2) \times\tfrac{1}{\beta}\Avc(1/2)
\]
such that $\{\Pi(\eta(t))\}_{t\in[0,1]}$ is $\C\left(\frac{2\beta}{4\beta-\gamma}\right)$.
The closed curve $\{\Pi(\eta(t))\}_{t\in[0,1]}$ strictly encloses $z$.
Since $\Avc(1/2) \times\Avc(1/2) \times\tfrac{1}{\beta}\Avc(1/2)$ is simply connected, we can continuously contract $\{\eta(t)\}_{t\in[0,1]}$ to a point in $\Avc(1/2) \times\Avc(1/2) \times\tfrac{1}{\beta}\Avc(1/2)$, and $\{\Pi(\eta(t))\}_{t\in[0,1]}$ continuously contracts to a point in $S_{\beta,\gamma}$.
However, this is not possible as $\{\Pi(\eta(t))\}_{t\in[0,1]}$ has a nonzero winding number around $z$ and $z\notin S_{\beta,\gamma}$.
We have a contradiction and we conclude $z\in S_{\beta,\gamma}$ and $\Avc\left(\frac{2\beta}{4\beta-\gamma}\right)\subseteq S_{\beta,\gamma}$.\qed

\section*{Acknowledgments}
This work was supported by NSF DMS-1720237 and ONR N000141712162.

\bibliography{my_srg,Master_Bibliography}

\end{document}

%% file: FancyMacros.tex
\usepackage{mathtools}
\usepackage[normalem]{ulem} 

\usepackage{amsmath}
\usepackage{amssymb}
\usepackage{mathtools}
\usepackage{mathrsfs}

\newcommand{\cut}[1]{{}}

\newcommand{\DeclareAutoPairedDelimiter}[3]{%
	\expandafter\DeclarePairedDelimiter\csname Auto\string#1\endcsname{#2}{#3}%
	\begingroup\edef\x{\endgroup
		\noexpand\DeclareRobustCommand{\noexpand#1}{%
			\expandafter\noexpand\csname Auto\string#1\endcsname*}}%
	\x}

\DeclareAutoPairedDelimiter{\p}{(}{)} 					
\DeclareAutoPairedDelimiter{\sp}{[}{]} 					
\DeclareAutoPairedDelimiter{\abs}{|}{|} 					
\DeclareAutoPairedDelimiter{\cp}{\{}{\}} 				
\DeclareAutoPairedDelimiter{\dotp}{\langle}{\rangle} 	
\DeclareAutoPairedDelimiter{\n}{\Vert}{\Vert} 			
\DeclareAutoPairedDelimiter{\cl}{\lceil}{\rceil}

\newcommand{\cA}{{\mathcal{A}}}

\newcommand{\cC}{{\mathcal{C}}}

\newcommand{\cG}{{\mathcal{G}}}

\newcommand{\cM}{{\mathcal{M}}}
\newcommand{\cN}{{\mathcal{N}}}

\newcommand{\cT}{{\mathcal{T}}}

